\documentclass[10pt]{article}
\usepackage{amsmath}
\usepackage{amsthm}
\usepackage{amsfonts}
\usepackage{amssymb}
\usepackage{latexsym} 

\usepackage[matrix,tips,graph,curve]{xy}



\makeatletter 
\@addtoreset{equation}{section}
\makeatother

\theoremstyle{plain}
\newtheorem{theorem}[equation]{Theorem}

\newtheorem{proposition}[equation]{Proposition}

\theoremstyle{definition}

\newtheorem{remark}[equation]{Remark}

\newtheorem{notation}[equation]{Notation}

\newcommand{\IC}{\mathbb{C}}

\newcommand{\IR}{\mathbb{R}}

\newcommand{\sign}{\rm sign}

\renewcommand\det{{\rm det\,}}

\def\d/{/\mspace{-6.0mu}/}

\newcommand{\db}{\bar{\partial}}

\begin{document}

\title{Fixed point thorems from a de Rham perspective}
\author{Mark Stern \footnote{Duke University, Department of Mathematics;  
 stern@math.duke.edu. The author was partially supported by NSF grant DMS-0504890.}
}

\date{\today}

\maketitle

\setcounter{section}{0}

\medskip
\section{Introduction}

Let $M$ be a smooth compact oriented Riemannian manifold of dimension $n$, and $f:M\rightarrow M$ a smooth map. Define 
the Lefschetz number
$$L(f) = \sum_{p=0}^n(-1)^pTrace(f^*|_{H^p(M)}).$$
The classical Lefschetz fixed point theorem states that if $f$ has isolated nondegenerate fixed points, then
$$L(f) = \sum_{f(b)=b}\sign\, \det(df_b-I).$$

Atiyah and Bott (\cite{AB1},\cite{AB2}) generalized this theorem to complexes of elliptic operators; we briefly recall (under mild restrictions) their theorem.
Let $E_0$, $E_1$, $\cdots, E_N$ be a sequence of smooth hermitian vector bundles over $M$, equipped with a sequence of first order differential operators $D_i:\Gamma(E_i)\rightarrow \Gamma(E_{i+1}).$  This sequence, denoted $\Gamma(E)$, is called an elliptic complex if for all $i$,
$$D_{i+1}D_i = 0,$$
  and 
$$D_i^*D_i + D_{i-1}D_{i-1}^*\mbox{     is elliptic}.$$
 Here we set $D_i = 0$ for $i\not\in [0,N-1].$ Set
$$H^p(\Gamma(E)) = Ker D_p/ImD_{p-1}.$$
Given a smooth map $f$ and smooth bundle homomorphisms $\phi_p:(f^*E)_p\rightarrow E_p$, we may define endomorphisms 
$T_p:\Gamma(E_p)\rightarrow \Gamma(E_p)$ by 
$$T_p s = \phi_p f^* s.$$
When 
\begin{equation}\label{intertwining}
D_pT_p = T_{p+1}D_p,
\end{equation}
 $T$ is called a geometric endomorphism of the complex $\Gamma(E)$. It induces endomorphisms $H^pT$ of $H^p(\Gamma(E)),$ and we define the Lefschetz number 
$$L(T) = \sum_{p=0}^N(-1)^pTrace\, H^pT.$$
Note that in the important special case where $N=1$ and $T$ is the identity map, $L(T) = index(D_0).$
The Atiyah-Bott theorem expresses the Lefschetz number in terms of fixed point data. 

\begin{theorem}(Atiyah-Bott, \cite{AB1}\cite{AB2}) Let $\Gamma(E)$ be an elliptic complex. Let $T$ be a geometric endomorphism of $\Gamma(E)$ associated to a pair $(f,\phi)$, with $f:M\rightarrow M$  a smooth map with isolated nondegenerate fixed points, and $\phi$ a smooth bundle homomorphisms $\phi_p:(f^*E)_p\rightarrow E_p$. Then 
$$L(T) = \sum_{f(b) = b}\frac{\sum_{p=0}^N(-1)^pTrace\,\phi_{p,b}}{|det(I-df_b)|}.$$
\end{theorem}

The proof of this theorem has gone through several incarnations. 
Set 
$$L_{ev}:=\oplus_{i}(D_{2i}^*D_{2i} + D_{2i-1}D_{2i-1}^*), \,\,\,\,\mbox{and}$$  
$$L_{od}:=\oplus_{i}(D_{2i+1}^*D_{2i+1} + D_{2i}D_{2i}^*).$$ The modern analytic proof of both the Atiyah-Bott theorem and the index theorem starts from the observation that the $\lambda$ eigenspaces of $L_{ev}$ and $L_{od}$ are isomorphic for $\lambda\not = 0$. The isomorphism is given by $\oplus_{i}\lambda^{-1/2}(D_{2i}+D_{2i-1}^*).$ Let $T_{ev}$ denote the restriction of $T$ to $\oplus_i\Gamma(E_{2i})$. Similarly define $T_{od}$. Then one can use this isomorphism and a Hodge isomorphism to show that 
\begin{equation}\label{intro1}
Trace\, T_{ev}e^{-tL_{ev}} - Trace\, T_{od}e^{-tL_{od}} = L(T),
\end{equation}
because the traces over the nonzero eigenspaces cancel identically. Then one can use elementary heat equation asymptotics to compute that 
\begin{equation}\label{intro2}
lim_{t\rightarrow 0}(Trace\,T_{ev}e^{-tL_{ev}} - Trace\, T_{od}e^{-tL_{od}}) = 
\sum_{f(b) = b}\frac{\sum_{p=0}^N(-1)^pTrace\,\phi_{p,b}}{|det(I-df_b)|}.
\end{equation}
The heat equation proof of the fixed point theorems was developed by numerous authors including 
Berline and Vergne \cite{BV}, Donnelly \cite{D1},\cite{DP}, Gilkey \cite{Gil2}, Kotake \cite{Ko}, Lafferty \cite{La}, Patodi \cite{P},\cite{DP}, and many others.  

In this note we wish to consider new extensions of the fixed point formulas. We pose the question: when can we extend the Lefschetz theorems to new classes of maps? More geometrically, passing from maps to their graphs, when can we extend the Lefschetz theorems to new classes of subspaces of $M\times M$?
The classical Lefschetz theorem applies to all maps, but only gives data associated to a single classical elliptic complex. If one further assumes that a map is an isometry, then one gains Lefschetz data associated to the signature complex. The greater the number of complexes for which one has Lefschetz data for a function $f$, the better one can analyze $f$. This paper was motivated by the observation that the index theoretic proof of the Lefschetz formulas often requires greater restrictions on the function $f$ (or more generally a correspondence $f$) than is actually necessary. For example, as we shall see in Theorem \ref{confsig}, a correspondence need only be locally conformal in order to be able to compute its Lefschetz number associated to the signature complex. Before stating further, more exotic new fixed point theorems, we explain the elementary ideas underlying the new results. 

Many older proofs of fixed point theorems (see \cite{T} and references therein), use an alternate analytic approach, closer in spirit to the original Lefschetz argument. Translating from the use of Green's operators in \cite{T} to the current preference for heat operators, this alternate proof exploits (\ref{intro1} and \ref{intro2}) but justifies (\ref{intro1}) slightly differently, from a de Rham perspective. The intertwining of the eigenspaces with nonzero eigenvalues is interpreted as follows. Let $\Gamma_f\subset M\times M$ denote the graph of $f$. Then $Trace\, T_{ev}e^{-tL_{ev}} - Trace\, T_{od}e^{-tL_{od}}$ can be computed as the integral  over $\Gamma_f$ of a differential form associated to the Schwartz kernels of the heat operators. We call this form {\em the index form}. Splitting the index form into a $t-$independent summand associated with projection onto the kernels of $L_{ev}$ and $L_{od}$ and a $t-$ dependent summand associated with all the nonzero eigenvalues, one finds that the $t-$dependent summand restricted to $\Gamma_f$ is exact. This exactness implies that the traces are $t-$independent, yielding (\ref{intro1}). 

In this note we seek larger classes of submanifolds (or more generally closed currents) $\Sigma$ in $M\times M$ on which the $t-$dependent summand of the index form becomes exact. The intertwining condition (\ref{intertwining}) is replaced by geometric conditions on $\Sigma$. For the index form corresponding to the Gauss-Bonnet index, the only condition on $\Sigma$ is that it be closed. For the Riemann-Roch index form of an $m$ complex dimensional compact manifold, $\Sigma$ can be a complex submanifold. In fact, we may consider more general currents than submanifolds. Then one requires the currents for the Riemann-Roch index form to be closed and have Hodge bitype $(m,m)$. In order to obtain Lefschetz formulas appropriate for families of complex maps, we may consider currents $T$ associated to a pair $(\Sigma,z)$ with $\Sigma$ an $m+p$ complex dimensional subvariety of $M\times M$ and $z$ a closed $(p,p)$ form. For the index form corresponding to the signature operator, $\Sigma$ is required to have the local structure of the graph of a conformal map. More generally, we may consider {\em extended conformal pairs} $(\Sigma,z)$ (see (\ref{cpair})) which provide Lefschetz formulas appropriate for families of conformal maps. 

In addition to these classical index forms, the search for geometries which are compatible with suitable index forms, where compatibility is understood in terms of the vanishing of the $t-$dependent summand,  leads to hybrid Lefschetz results. For example, coisotropic geometry suggests the introduction of a new index form corresponding to $Tr J\ast e^{-t\Delta}$, where $J$ is a an almost complex structure operator and $\ast$ is the Hodge star operator. Applied to $\Sigma = V_1\times V_2$, with $V_i$ coisotropic, we obtain intersection formulas, which, however, reduce to formulas obtainable from the Gauss-Bonnet index form (and are therefore effectively classical). 

In a different direction, we examine the elimination of the $t-$ dependent terms in the index form, via averaging rather than geometry. In the examples we consider, this requires some homogeneous structure. For example, on an abelian variety, $A$, we can integrate the Riemann-Roch index form over  a family 
$\{\Sigma_y = V_1\times V_{2y}\}_{y\in A}$, with $V_1$ and $V_{2y}$ special Lagrangian varieties. We obtain formulas relating cohomological data to average intersection data for the family. We obtain similar results for families of Lagrangian varieties in compact hermitian symmetric spaces. 

We doubt that we have come close to exhausting the possible applications of these elementary ideas. We have worked in the combinatorially trivial regime of transverse intersections of $\Sigma$ with the diagonal of $M\times M$. This corresponds to considering only isolated nondegenerate  fixed points. There is no apparent obstruction to using the calculus of Clifford modules and Mehler's formula (as in for example \cite{BGV},\cite{La}) to treat the more general case, but we do not pursue that direction here. 
 
In the following, we will first recall this de Rham perspective for the classical de Rham and Dolbeault
complexes. We will then show how to apply it to obtain new fixed point theorems for conformal relations, coisotropic intersections, and average special Lagrangian intersections. We end with an appendix in which we adapt standard heat equation asymptotics to our context. 

\section{Gauss Bonnet}
In this section we illustrate the de Rham perspective on index theory in the simplest possible case, the Gauss-Bonnet theorem. All the results in this section are well known. Let $M^n$ be a compact oriented riemannian manifold. 
Let $\Delta_p$ denote the Laplace Beltrami operator on the space $A^p$ of $p-$forms. Let $e^p_t(x,y)$ denote the Schwartz kernel for $e^{-t\Delta_p}.$  We define the Schwartz $n$ form $e^p_t(x,y)$ by requiring for every $p-$form $f$,
$$e^{-t\Delta_p}f(x) = \int_{\{x\}\times M}e^p_t(x,y)\wedge f(y).$$
Here we adopt the convention that when treating forms on $M\times M$, $f(x)$ and $f(y)$ are shorthand for $\pi_1^*f$ and $\pi_2^*f$ respectively, where $\pi_1$ and $\pi_2$ denote the projections on the first and second factors of $M\times M$ respectively. Throughout this note the Hodge star operator, unless otherwise subscripted, will denote the Hodge star operator for $M$ rather than for $M\times M$. Thus on $M\times M$, $\ast f(y)$ denotes $\pi_2^*(\ast f)$ and similarly for $\ast f(x)$. 

Expanding $e^p_t$ with respect to an orthonormal basis of eigenforms $\{\phi^p_{\lambda}\}$ of $\Delta$ we have 
$$e^p_t(x,y) = (-1)^{p(n-p)}\sum_{\lambda}e^{-t\lambda}\phi^p_{\lambda}(x)\wedge \ast\phi^p_{\lambda}(y).$$
Then 
$$Tr\, e^{-t\Delta_p} = \int_{\delta} (-1)^{p(n-p)}e^p_t,$$
where $\delta $ denotes the diagonal in $M\times M$. The Euler characteristic of $M$ is given  by 
$$\sum_p(-1)^pTr\, e^{-t\Delta_p} = \int_{\delta} \sum_p(-1)^{p(n-p+1)}e^p_t.$$
We call $\sum_p(-1)^{p(n-p+1)}e^p_t$ the {\em index form} associated to the Gauss Bonnet theorem. 

Refining the eigen-expansion further in terms of harmonic, closed and coclosed forms, we have 
$$\sum_p(-1)^{p(n-p+1)}e^p_t = 
(-1)^{p}\sum_{i,p}h^p_{i}(x)\wedge \ast h^p_i(y) +
(-1)^{p}\sum_{\lambda,p}e^{-t\lambda}\lambda^{-1}d\psi^{p-1}_{\lambda}(x)\wedge \ast d\psi^{p-1}_{\lambda}(y)$$
$$
+ (-1)^{p}\sum_{\lambda,p}e^{-t\lambda}\psi^p_{\lambda}(x)\wedge \ast\psi^p_{\lambda}(y),$$
where we take $\{\psi^p_{\lambda}\}$ to be an orthonormal coclosed eigenbasis and $\{h_i\}$ an orthonormal basis of harmonic forms. 
Using the standard relation on $k-$ forms: 
$$d^*f = (-1)^{k}\ast^{-1} d\ast,$$
we regroup this expansion as 
\begin{equation}\label{GB1}
\sum_p(-1)^{p(n-p+1)}e^p_t = 
\sum_{i,p}(-1)^{p}h^p_{i}(x)\wedge \ast h^p_i(y) +
d(\sum_{\lambda,p}(-1)^{p}e^{-t\lambda}\lambda^{-1}\psi^{p-1}_{\lambda}(x)\wedge \ast d\psi^{p-1}_{\lambda}(y)).
\end{equation}
Thus the $t-$dependent summand of the index form is exact. This reduces the Gauss Bonnet and Lefschetz fixed point theorems to standard heat equation asymptotics, which we now recall. (See for example \cite{BGV} or \cite{Gil}). Define the $n-$ form $\nu$ on the diagonal, $\delta$, to be the volume form for the normal bundle $N_{\delta}$ to $\delta$. In particular, if $\mu^j$ is a local oriented orthonormal coframe on $M$ near $b$, then near $(b,b)\in\delta$ 
\begin{equation}\label{nu}
 \nu:=2^{-n/2}(\pi_1^*\mu^1-\pi_2^*\mu^1)\wedge\cdots\wedge (\pi_1^*\mu^n-\pi_2^*\mu^n).
\end{equation}
Given a submanifold $W^{n+c}\subset M\times M$ intersecting the diagonal transversely in a submanifold $S$, we let $N_S^W$ denote the normal bundle to $S$ in $W$. The assumption of transversality implies that dim $N_S^W = n$ and that the projection 
$\Pi_W(b):N_{S,(b,b)}^W\rightarrow N_{\delta}$ is an isomorphism. With this notation, heat equation asymptotics (see \ref{append}) give for $\Sigma$ an $n$ dimensional submanifold intersecting $\delta$ transversely, 
\begin{equation}\label{GBa}
lim_{t\rightarrow 0}\int_{\Sigma}\sum_p(-1)^{p(n-p+1)}e_t^p = \sum_{(b,b)\in V\cap\delta}\frac{\langle2^{n/2}\nu,dV_{\Sigma}\rangle}{det^{1/2}(\Pi_W^*\Pi_W(b))}.
\end{equation}
Here $dV_{\Sigma}$ denotes the volume form of the submanifold $\Sigma$,  identified via the metric with an element of $\bigwedge^nT^*(M\times M)_{|\delta}$. Setting 
$$\nu_{GB}(T_{(b,b)}\Sigma ) = \frac{\langle 2^{n/2}\nu,dV_{\Sigma}\rangle}{det^{1/2}(\Pi_W^*\Pi_W(b))},$$
we obtain 
\begin{theorem}\label{GBL}(Lefschetz)
Let $M^n$ be a smooth compact oriented Riemannian manifold.
Let $\Sigma $ be a smooth n dimensional submanifold of $ M\times M$. Assume that $\Sigma$ intersects $\delta$ transversely. Then 
\begin{equation}\label{GBLe}
\int_{\Sigma}\sum_{i,p}(-1)^{p}h^p_{i}(x)\wedge \ast h^p_i(y) = \sum_{(b,b)\in \Sigma\cap \delta}\nu_{GB}(T_{(b,b)}\Sigma ).
\end{equation}
\end{theorem}
The proof of this result consists of integrating equation (\ref{GB1}) over $\Sigma$. The $t$-dependent terms are exact, giving us the left hand side of (\ref{GBLe}). On the other hand, taking the limit as $t\rightarrow 0$ and applying (\ref{GBa}) gives the righthand side.  In  (\ref{example})  we show that when $\Sigma$ is locally the graph of a function $f$, then as expected,
$$\nu_{GB}(T_{(b,b)}\Sigma ) = sign\,det(I-df).$$

\section{Riemann Roch}Next we turn to the Dolbeault complex. (See \cite{TT} for an extensive treatment of the Lefschetz theorem for this complex.)  Let $M$ be a compact complex manifold of complex dimension $m$. 
Define $e_t^{0,q}$ to be the Schwarz kernel double form for $e^{-t\Box_q}$, where $\Box_q$ denotes the Dolbeault Laplacian on $(0,q)$ forms. Thus, in the kahler case it is (up to rescaling $t$) the summand of $e_t^q$ obtained by projecting $e^{-t\Delta_q}$ onto $(0,q)$ forms. In this case $n = dim_{\IR}M$ is even, and 
$$e^{0,q}_t(x,y) = (-1)^{q}\sum_{\lambda}e^{-t\lambda}\phi^{0,q}_{\lambda}(x)\wedge \bar\ast\phi^{0,q}_{\lambda}(y).$$
We call $\sum_q e^{0,q}_t(x,y)$ the {\em index form} for the Dolbeault complex. We similarly refer to the index form associated to other index problems. Let $A^{p,q}$ denote the space of differential forms of bitype $(p,q)$. 

We first examine whether the cohomology class of the Dolbeault index form is $t$ independent. We expand it as 
$$\sum_q e^{0,q}_t(x,y) =  
\sum_q(-1)^{q}\sum_{i}h_i^{0,q}(x)\wedge \bar\ast h_i^{0,q}(y) +
\sum_q(-1)^{q}\sum_{\lambda}e^{-t\lambda}\lambda^{-1}\db b^{0,q-1}_{\lambda}(x)\wedge \bar\ast\db b^{0,q-1}_{\lambda}(y)$$
$$ + \sum_q(-1)^{q}\sum_{\lambda}e^{-t\lambda}b^{0,q}_{\lambda}(x)\wedge \bar\ast b^{0,q}_{\lambda}(y),$$
where the $b_{\lambda}$ are an eigenbasis for the co-$\db$-exact forms. We rewrite this as 

$$ \sum_q(-1)^{q}\sum_{i}h_i^{0,q}(x)\wedge \bar\ast h_i^{0,q}(y) +
\sum_q(-1)^{q}d(\sum_{\lambda}e^{-t\lambda}\lambda^{-1} b^{0,q-1}_{\lambda}(x)\wedge \bar\ast\db b^{0,q-1}_{\lambda}(y))$$
$$-\sum_q(-1)^{q}\sum_{\lambda}e^{-t\lambda}\lambda^{-1} \partial b^{0,q-1}_{\lambda}(x)\wedge \bar\ast\db b^{0,q-1}_{\lambda}(y).$$
Thus we see that the Dolbeault index form is the sum of a t independent form, an exact form, and a form of type $(m+1,m-1)$. Hence the cohomology class of this form is apparently {\em not} $t-$independent. Let $\Sigma$ be a middle dimensional submanifold of $M\times M$. Let $T_{\Sigma}$ denote the corresponding current. If $T_{\Sigma}$ is of type $(m,m)$ then the integral of the index form over $\Sigma$  {\em is} independent of $t$. Every m dimensional complex subvariety of $M\times M$ determines a current of this type. A theorem of King \cite[Theorem 5.2.1]{K} and its generalization by Harvey and Shiffman \cite{HS} imply that there are essentially no additional geometric examples.  Choosing $\Sigma$ to be a correspondence, we obtain the holomorphic Lefschetz theorem for correspondences. From the appendix, we obtain the following expression for the heat equation asymptotics. 
\begin{equation}\label{nRR}
lim_{t\rightarrow 0}\int_{\Sigma}\sum_qe_t^{0,q} = \sum_{(b,b)\in \Sigma\cap\delta}\frac{\langle 2^{m/2}\nu_2,dV_W\rangle}{det^{1/2}(\Pi_{\Sigma}^*\Pi_{\Sigma})},
\end{equation}
where in a local unitary frame, $\{\eta^i\}_i$, of the holomorphic cotangent bundle 
\begin{equation}\label{nu2}
2^{m/2}\nu_2 := (-1)^{m(m+1)/2}(\pi_1^*\bar \eta^1-\pi_2^*\bar\eta^1)\wedge\cdots\wedge (\pi_1^*\bar \eta^m-\pi_2^*\bar\eta^m)\wedge \pi_2^*\eta^1\wedge\cdots\wedge \pi_2^*\eta^m.
\end{equation}
Set 
\begin{equation}\label{nRR2}
\nu_{RR}(T_{(b,b)}\Sigma) := \frac{\langle 2^{m/2}\nu_2,dV_W\rangle}{det^{1/2}(\Pi_{\Sigma}^*\Pi_{\Sigma})}.
\end{equation}
 When $\Sigma$ is the graph of a holomorphic function $f$, then a computation like that of (\ref{example}) gives
$$\nu_{RR}(T_{(b,b)}\Sigma) = \frac{1}{det_{\IC}(I-df_b)}.$$
Combining the $t$ independence of the integral of the Riemann-Roch index form on $\Sigma$ with the heat equation asymptotics, we obtain the following holomorphic Lefschetz fixed point theorem.
\begin{theorem}
Let $M$ be a compact Kahler manifold.
Let $\Sigma\subset M\times M$ be a holomorphic correspondence. Let 
$\{h_i^{0,p}\}_i$ be an orthonormal basis of harmonic $(0,p)$ forms on $M$. 
Then 
$$\int_{\Sigma}\sum_p(-1)^{p}\sum_{i}h_i^{0,p}(x)\wedge \bar\ast h_i^{0,p}(y) = \sum_{(b,b)\in \Sigma\cap\delta}\nu_{RR}(T_{bb}\Sigma).$$
\end{theorem}
When $\Sigma$ is the graph of a holomorphic function $f$, this reduces to the usual holomorphic Lefschetz theorem. 

At the opposite extreme, we may take $\Sigma = V\times W$, with $V$ and $W$ smooth subvarieties of $M$ of complementary dimension. Then $V\times W$ is a middle dimensional subvariety of type $(m,m)$ of $M\times M$. This gives a trivial result since the index form identically vanishes on this variety. So, to recapture standard intersection results, we simply replace $(0,q)$ forms by $(p,q)$ forms, with $p= dim_C V$, and then repeat the preceding computations.  
\subsection{Excess dimensions}\label{excesssub}
Suppose that $z$ is a $\db$ closed $(p,p)$ form on $M\times M$. Then $z\wedge \sum_qe^{0,q}_t$ once again decomposes into a $t-$ independent summand $z\wedge \sum_q(-1)^q\sum_ih_i^{0,q}(x)\wedge \bar\ast h_i^{0,q}(y),$ an exact summand, and a summand with no $(m+p,m+p)$ component. So, choosing an $m+p$ dimensional complex subvariety $W$ of $M\times M$, we obtain once again an equality between cohomological data and fixed point data. Let $S=W\cap\delta$. Restricted to (as opposed to pulled back to) $\delta$, we can use the metric to factor orthogonally  
$$dV_W = \nu_S^W\wedge dV_S,$$
with $\nu^W_S$ a volume form for the fiber of the normal bundle $N_S^W$ of $S$ in $W$. 
We set 
$$\nu_{RR}(z,T_{(b,b)}W) =  \frac{\langle z(b,b)\wedge 2^{m/2}\nu_2,dV_W(b)\rangle }{det^{1/2}(\Pi^*_{W,b}\Pi_{W,b})}=  \frac{\langle z(b,b)\wedge 2^{m/2}\nu_2,\nu^W_S\wedge dV_S\rangle}{det^{1/2}(\Pi^*_{W,b}\Pi_{W,b})}.$$
Only the  $(p,a)(0,p-a)$ component of $z$ contributes to this expression. Here the multigrading is the refined Hodge grading associated with the product structure. In particular, a form $\sum_{|I|=a,|J|=b,|K|=c,|L|=d}f_{IJKL}\pi_1^*\omega^I\wedge \pi_1^*\bar \omega^J\wedge \pi_2^*\omega^K\wedge \pi_2^*\bar \omega^L$ has type $(a,b)(c,d)$. 
\begin{notation}Given a $(p,p)$ form $\phi$, let  $\phi_0$ denote the sum over $a$ of its $(p,a)(0,p-a)$ components. 
\end{notation}
Then 
$$\langle z(b,b)\wedge 2^{m/2}\nu_2,\nu^W_S\wedge dV_S\rangle = 2^p\langle z_0(b,b), dV_S\rangle \langle2^{m/2}\nu_2,\nu^W_S\rangle,$$
and 
$$\nu_{RR}(z,T_{(b,b)}W) = 2^p\frac{\langle 2^{m/2}\nu_2,\nu^W_S\rangle\langle z_0,dV_S\rangle}{det^{1/2}(\Pi^*_{W,b}\Pi_{W,b})}.$$
Then our arguments (see Example 3 in the appendix) yield in this case:
\begin{theorem}
Let $W$ intersect $\delta$ transversely in a smooth $p$ dimensional complex submanifold, $S$. Then 
if $z$ is of type $(p,p)$,
$$\int_W z(x,y)\wedge \sum_q(-1)^q\sum_ih_i^{0,q}(x)\wedge \bar\ast h_i^{0,q}(y) = \int_S 
\frac{2^p\langle 2^{m/2}\nu_2,\nu^W_S\rangle z_0}{det^{1/2}(\Pi^*_{W,b}\Pi_{W,b})}.$$
\end{theorem}
\begin{remark}
The reader familiar with the Atiyah-Segal-Singer Lefschetz theorems for fixed submanifolds of dimension greater than zero (\cite{ASe},\cite{ASi}) may be surprised that no curvature data explicitly enters into the above fixed point formula. As we show in the appendix, it is not the dimension of $S$ in $\delta$ but the codimension of $S$ in $W$ which is germane for determining the simplicity of the fixed point data. If the codimension is equal to the dimension of $M$, then no difficult combinatorics or curvature computations are required. 
\end{remark}
Suppose now that $M$ is Kahler and that $z$ is $d-$closed and dual to a cycle $V$ in $M\times M$. Suppose further that $V$ is a complex subvariety. Let $T_V$ be the current corresponding to integrating over $V$. Then 
$$z-T_V = \partial\db b$$
 for some $(m-p-1,m-p-1)$ current $b$. Then 
   
$$\int_{W} z\wedge \sum_q(-1)^{q}\sum_{i}h_{i}^{0,q}(x)\wedge \bar\ast h_{i}^{0,q}(y)
 = \int_{V\cap W} \sum_q(-1)^{q}\sum_{i}h_{i}^{0,q}(x)\wedge \bar\ast h_{i}^{0,q}(y).$$

On the other hand, for $V$ holomorphic,
$$\int_{V\cap W} \sum_q e^{0,q}_t(x,y) = \int_{V\cap W} \sum_q(-1)^{q}\sum_{i}h_{i}^{0,q}(x)\wedge \bar\ast h_{i}^{0,q}(y).$$

Hence if $V\cap W = Y$ is transverse to the diagonal, then we have 

\begin{equation}\label{equalloc}
\int_S \nu_{RR}(z,T_{(b,b)}W)dV_S = \sum_{(b,b)\in Y\cap \delta}\nu_{RR}(T_{(b,b)}(Y)).
\end{equation}
We can reexpress this as 

\begin{equation}\label{equalloc2}
\int_S \frac{2^p\langle 2^{m/2}\nu_2,\nu^W_S\rangle z_0}{det^{1/2}(\Pi^*_{W,b}\Pi_{W,b})} = \sum_{(b,b)\in V\cap S}\frac{\langle 2^{m/2}\nu_2,\nu^{V\cap W}_{(b,b)}\rangle }{det^{1/2}(\Pi^*_{V\cap W,b}\Pi_{V\cap W,b})}.
\end{equation}

For example, if $W = M\times M$, ($M$ connected) this reduces to  
\begin{equation}\label{equalloc4}
\int_\delta z_0  = \int_{M\times M} z\wedge \sum_{q,i}(-1)^{q}h_{i}^{0,q}(x)\wedge \bar\ast h_{i}^{0,q}(y) = \sum_{(b,b)\in V\cap \delta}\frac{\langle 2^{m/2}\nu_2,\nu^{V}_{(b,b)}\rangle }{det^{1/2}(\Pi^*_{V,b}\Pi_{V,b})}.
\end{equation}

Observe that equality (\ref{equalloc2}) depends on both the hypothesis that $z$ is type $(p,p)$ and the additional assumption that $V$ is complex. 
The apparent difference in the 2 sides of the formula is interesting in light of the relation between the two hypotheses and the Hodge conjecture.  In order 
to extract more information from such relations, it would be useful to have an infinite family of such formulas for fixed $V$, $z$, and $W$. Hence we next consider holomorphic Lefschetz formulas for the Dolbeault complex with coefficients in a holomorphic vector bundle. 

\subsection{Holomorphic Coefficients}

Consider  now the $\db$ operator with coefficients in a hermitian holomorphic vector bundle, $E$. Then the  heat kernel becomes 
$$e^{0,q}_t(x,y) = (-1)^{q}\sum_{\lambda}e^{-t\lambda}\phi^{0,q}_{\lambda,a}(x)\wedge \bar\ast\phi^{0,q}_{\lambda,b}(y)s^a(x)(\cdot,s^b(y)),$$
where $\{s^a\}_a$ is a local holomorphic frame and $\phi_{\lambda} = \phi_{\lambda,a}s^a$ in this frame. In order to obtain a scalar valued differential form, we need additional data in the form of a holomorphic section $Q$ of 
$Hom(\pi_1^*E\otimes \pi_2^*E^*,\IC)$ over a subvariety $W$. (In the case of tensor bundles such sections arise as the pullback to $W$ under the Gauss map of a global section of 
$Hom(\pi_G^*(\pi_1^*E\otimes \pi_2^*E^*),\IC),$ where 
$\pi_G: G_n(T(M\times M))\rightarrow M\times M$ denotes the projection from the grassmanian of $n$ planes in the tangent bundle.)  Set $k^{ab}(x) = \langle s^a(x),s^b(x)\rangle$. Let $\{w_a\}_a$ denote the holomorphic coframe dual to $\{s^a\}_a$. Let $Q^a_c = Q(s^a\otimes w_c).$ This function is a local holomorphic function on $W$. Then we can generalize our earlier expansions, writing 
$$\sum_q Q(e^{0,q}_t(x,y)) =  
\sum_q(-1)^{q}\sum_{i}Q^a_ch_{i,a}^{0,q}(x)\wedge \bar\ast h_{i,b}^{0,q}(y)k^{bc}$$
$$ +
\sum_{q,\lambda}(-1)^{q}e^{-t\lambda}\lambda^{-1}Q^a_c\db b^{0,q-1}_{\lambda,a}(x)\wedge \bar\ast\db b^{0,q-1}_{\lambda,b}(y)k^{bc} + \sum_{q,\lambda}(-1)^{q}e^{-t\lambda}Q^a_cb^{0,q}_{\lambda,a}(x)\wedge \bar\ast b^{0,q}_{\lambda,b}(y)k^{bc},$$
where the $b_{\lambda}$ are again an eigenbasis for the co $\db$ exact forms. We rewrite this (on $W$) as 
$$\sum_{q,i}(-1)^{q}Q^a_ch_{i,a}^{0,q}(x)\wedge \bar\ast h_{i,b}^{0,q}(y)k^{bc} +
d\sum_{q,\lambda}(-1)^{q}e^{-t\lambda}\lambda^{-1} Q^a_cb^{0,q-1}_{\lambda,a}(x)\wedge \bar\ast\db b^{0,q-1}_{\lambda,b}(y)k^{bc}$$
$$-\sum_{q,\lambda}(-1)^{q}e^{-t\lambda}\lambda^{-1} \partial (Q^a_cb^{0,q-1}_{\lambda,a}(x))\wedge \bar\ast\db b^{0,q-1}_{\lambda,b}(y)k^{bc}$$
$$-\sum_{q,\lambda}(-1)^{q}e^{-t\lambda}\lambda^{-1}\db Q^a_c\wedge b^{0,q-1}_{\lambda,a}(x)\wedge  \bar\ast \db b^{0,q-1}_{\lambda,b}(y)k^{bc}.$$
 
The assumption that $Q$ is holomorphic on $W$ eliminates the last term, and we are left again with the index form restricted to $W$ as the sum of a t independent form, an exact form, and a form of type $(m+1,m-1)$, where $m$ denotes the complex dimension. Consequently, the integral of this form over the subvariety $W$ of $M\times M$ is independent of $t$.   

More generally, we can consider a $\db$ closed $Q\in A^{r,r}(W,Hom(\pi_1^*E\otimes \pi_2^*E^*,\IC)),$ for any $q$ and repeat the preceding computation to obtain the following generalization of the holomorphic Lefschetz theorem. 
\begin{theorem}
Let $W$ be a $r+m$ dimensional subvariety of $M\times M$ intersecting $\delta$ transversely in a submanifold $S$. Let $Q\in A^{r,r}(W,Hom(\pi_1^*E\otimes \pi_2^*E^*,\IC))$ be $\db$ closed. Then 
$$\int_W\sum_q(-1)^{q}\sum_{i}Q^a_c\wedge h_{i,a}^{0,q}(x)\wedge \bar\ast h_{i,b}^{0,q}(y)k^{bc} = \int_S\nu_{RR}(trQ,T_{(b,b)}W)dV_S .$$
\end{theorem}

Specialize this theorem now to the case where $Q= z\wedge Q_1$, with $z$ a harmonic $(p,p)$ form. Suppose further that $z$ is dual to a complex submanifold $V\subset M\times M$ intersecting the diagonal transversely. Then arguing as in the preceding section and retaining its notation we have 

\begin{theorem}
Let $W$ be an $r+m$ dimensional subvariety of $M\times M$ intersecting $\delta$ transversely in a submanifold $S$. Let $Q_1\in A^{r,r}(W,Hom(\pi_1^*E\otimes \pi_2^*E^*,\IC))$ be $\db$ closed. Then 
$$\int_W\sum_q(-1)^{q}\sum_{i}z\wedge Q^a_{1c}\wedge h_{i,a}^{0,q}(x)\wedge \bar\ast h_{i,b}^{0,q}(y)k^{bc} = \int_S\nu_{RR}(z\wedge trQ_1,T_{(b,b)}W)dV_S$$
$$=   \int_{S\cap V}\nu_{RR}(trQ_1,T_{(b,b)}(W\cap V))dV_{S\cap V}.$$
\end{theorem}

\section{Signature}
Let $M$ be a $4k$ dimensional compact oriented Riemannian manifold. Let $\tau$ denote the involution of the exterior forms defined by Clifford multiplication by the volume form. In particular, if $\{\omega^i\}_i$ is a local orthonormal oriented coframe, and if we let $e(\omega^i)$ denote exterior multiplication on the left by $\omega^i$, then 
$$\tau = (e(\omega^1)-e^*(\omega^1))\cdots (e(\omega^{4k})-e^*(\omega^{4k})).$$
Recall that the topological signature of $M$ is equal, by the Hodge theorem, to the trace of $\tau$ restricted to the space of harmonic forms of $M$. The index form for the signature theorem is given by 

$$s_t = \sum_p (-1)^{p}\tau_x e^p_t(x,y).$$
Here $\tau_x$ denotes clifford multiplication by $\pi_1^*dvol$.
The de Rham explanation for the $t-$independence of the cohomology class of $s_t$ pulled back to the diagonal is somewhat different from the cases of the Gauss-Bonnet and Riemann Roch complexes, and perhaps should be called the Chern-Hirzebruch-Serre explanation (see \cite[Lemma 4]{CHS}). 
As usual, we expand the index form in terms of eigenfunctions, obtaining 
$$s_t =  \sum_{p,i}\tau h_{i}^p(x)\wedge \ast h_{i}^p(y)
+ \sum_{p,\lambda}e^{-t\lambda}\lambda^{-1}\tau d\psi_{\lambda}^{p-1}(x)\wedge \ast d\psi_{\lambda}^{p-1}(y)$$
$$
+ \sum_{p,\lambda}e^{-t\lambda}\lambda^{-1}\tau \ast d\psi_{\lambda}^{n-p}(x)\wedge \ast^2d\psi_{\lambda}^{n-p}(y).$$

The last sum is clearly exact. The second sum is exact on the diagonal because for all $2k$ forms $f$ and $F$, $\ast f(x)\wedge \ast F(y) - f(x)\wedge F(y)$ vanishes when pulled back to the diagonal. Forms of degree unequal to $2k$ do not contribute to the computation for dimension reasons. Therefore, the class of $s_t$ pulled back to the diagonal is $t-$ independent. Integrating $s_t$ on the diagonal then gives the signature theorem. In order to find new Lefschetz type theorems, we require a class of submanifolds of $M\times M$ on which 
$$B_t:= \sum_{p,\lambda}e^{-t\lambda}\lambda^{-1}\tau d\psi_{\lambda}^{p-1}(x)\wedge \ast d\psi_{\lambda}^{p-1}(y)$$ integrates to zero. Let $V$ be a closed $n-$dimensional submanifold of $M\times M$. A sufficient criterion to ensure that $B_t$ integrates to zero along $V$ is suggested by the diagonal. Call $V$ {\em self dual in middle dimension} if it satisfies the condition that 
for all $2k$ forms $f$ and $F$, $\ast f(x)\wedge \ast F(y) - f(x)\wedge F(y)$ vanishes when pulled back to $V$. We remark that if we did not restrict the degree of the forms required to satisfy this condition, then the current $T_V$ associated to $V$ would satisfy $dT_V = d^*T_V = 0$, which has no submanifold solutions, by elliptic regularity.   
  
Let $i_V:V\rightarrow M\times M$ denote the inclusion map. Suppose $(x,y)\in V$, and let $\{\omega^i(x)\}_i$ and $\{\omega^i(y)\}_i$ be any orthonormal frames for $\pi_1^*T^*M$ and $\pi_2^*T^*_YM$.  In order to satisfy the self duality condition, we need 
$$i_V^*(\omega^I(x)\wedge \omega^J(y)) = i_V^*(\ast \omega^I(x)\wedge \ast \omega^J(y)).$$
At a point $(x,y)\in V$ where the maps $i_{V}^*\pi_1^* :T_x^*M\rightarrow T_{(x,y)}^*V$ and $i_{V}^*\pi_2^*:T_y^*M\rightarrow T_{(x,y)}^*V$ are bijective, let 
$f = (i_{V}^*\pi_2^*)^{-1}i_{V}^*\pi_1^*.$ Choose now $\{\omega^i(x)\}_i$ to be an orthonormal $f^*f$  eigen basis for $T_x^*M$, with eigenvalues $\lambda^i$. For multiindices $I$, let $\lambda^I = \prod_{i\in I}\lambda^i$. Choose $\omega^i(y) = \dfrac{f\omega^i(x)}{\sqrt{\lambda^i}}$. Then 
$$i_V^*(\omega^I(x)\wedge \omega^J(y)) = \dfrac{1}{\sqrt{\lambda^J}}i_V^* \omega^I(x)\wedge i_V^*\pi_2^* f\omega^J(x) 
= \dfrac{1}{\sqrt{\lambda^J}}i_V^* (\omega^I(x)\wedge \omega^J(x)).$$
Hence the condition that $V$ be self dual in middle dimension becomes 
$$\lambda^I = \lambda^{I^c},$$
for all multi indices $I$ with $2|I| = dim V.$
This is equivalent to $\lambda^i = \lambda$, for some $\lambda$ independent of $i$.
In particular, this condition is satisfied if $f$ is locally the graph of a conformal morphism. We call such $V$ a {\em conformal correspondence}. 
 Computing the local heat asymptotics from (\ref{append}) yields 

\begin{equation}
\lim_{t\rightarrow 0}\int_Vs_t = \sum_{(b,b)\in V\cap\delta}\frac{\sum_{|I|=2k}\langle \pi_1^*\omega^I\wedge \pi_2^*\omega^I,dV_V\rangle}{det^{1/2}(2\Pi_V^*\Pi_V)}.
\end{equation}
We set 
$$\nu_{sig}(T_{(b,b)}V) := \frac{\sum_{|I|=2k}\langle \pi_1^*\omega^I\wedge \pi_2^*\omega^I,dV_V\rangle}{det^{1/2}(2\Pi_V^*\Pi_V)}.
$$

Then the $t-$independence of the integral of $s_t$ pulled back to $V$ gives the following. 
\begin{theorem}\label{confsig}
Let $V$ be a conformal correspondence in $M^{2m}\times M^{2m}$. Then 
$$\int_Vh_i(x)\wedge h_i(y) = \sum_{b\in V\cap \delta}\nu_{sig}(T_{(b,b)}V),$$
\end{theorem}

If $V$ is locally the graph of a conformal map $f$, we can write, for $(b,b)\in V\cap \delta$,
$df_b = \mu k$, with $\mu$ a positive scalar and $k$ an isometry,  with rotation angles $\theta_j$. Then block diagonalizing, we find 
$$\nu_{sig}(T_{(b,b)}V) = i^{-m}\frac{\mu^m|det(I-k)|}{|det(I-\mu k)|}\prod_{j=1}^mcot(\theta_j/2).$$

The extension of this theorem to the signature operator with coefficients in a flat bundle $E$ is more complex than in the Riemann-Roch case. If $Q$ is a flat section of $\pi_1^*E\otimes \pi_2^*E^*$ over a conformal correspondence $V$, then we have for $E$ unitary, by the preceding argument:
$$\int_VQ^{ab}h_{i,a}(x)\wedge h_{i,b}(y) = \sum_{b\in V_{\delta}}trQ\,\nu_{sig}(T_{(b,b)}V).$$
The unitary assumption can be removed at the cost of introducing additional notation; for simplicity, we will not address the nonunitary case. (See \cite{Lu}.)
The extension to closed differential forms, $Q$, of degree $p>0$, with coefficients in $\pi_1^*E\otimes \pi_2^*E^*$ and submanifolds $W$ of $M\times M$ of dimension $n+p$ requires a comaptibility condition between $W$ and $Q$ extending the conformal correspondence structure. In simplest form the compatibility condition becomes 
\begin{equation}\label{cpair}
Q\wedge i_W^*[f(x)\wedge F(y) - \ast f(x)\wedge \ast F(y)] = 0.
\end{equation}
We call a pair $(W,Q)$ satisfying $(\ref{cpair})$ an {\em extended conformal pair}. Such pairs can be constructed, for example,  when $W$ locally has the form of a family of conformal maps of $M$ over a base $B$ and $Q$ locally has the form of a pullback of a closed form on $B$. 
Define for a $p-$ form $z$,
$$\nu_{sig}(z,T_{(b,b)}W) := \frac{\sum_{|I|=2k}\langle z\wedge \pi_1^*\omega^I\wedge \pi_2^*\omega^I,dV_W\rangle}{det^{1/2}(2\Pi_W^*\Pi_W)}.
$$
We obtain the following theorem. 

\begin{theorem}\label{confsig}
Let $E$ be a unitary flat vector bundle over $M$. Let $Q$ be a closed $\pi_1^*E\otimes \pi_2^*E^*$ valued $p-$ form over a submanifold $W$ of $M\times M$. Assume that $(W,Q)$ is an extended conformal pair and $W$ intersects $\delta$ transversely. Then  
$$\int_WQ^{ab}\wedge h_{i,a}(x)\wedge h_{i,b}(y) = \int_{W\cap \delta}\nu_{sig}(trQ,T_{(b,b)}W)dV_{W\cap\delta}.$$

\end{theorem}

The extension to excess dimensions is, of course, easier for the Gauss-Bonnet index form. We have the following theorem. 
\begin{theorem}\label{fatgb}Let $M^n$ be a compact Riemannian manifold. Let $E$ be a unitary flat vector bundle over $M$. Let $\Sigma$ be an $n+p$ dimensional submanifold of $M\times M$, intersecting the diagonal transversely. 
Let $Q\in A^p(\Sigma,\pi_1^*E\otimes\pi_2^*E)$ be $d-$closed. Then 
$$\int_{\Sigma}Q^{ab}\wedge (-1)^{p}\sum_{i,p}h^p_{i,a}(x)\wedge \ast h^p_{i,b}(y)  = 
\int_{\Sigma\cap \delta }\frac{\langle tr\,Q\wedge \nu,dV_{\Sigma}\rangle dV_{\Sigma\cap \delta}}{det^{1/2}(\Pi^*_{\Sigma}\Pi_{\Sigma})}.$$
\end{theorem}

As we see in Example 2 in the appendix,
$$\int_{\Sigma\cap \delta }\frac{\langle tr\,Q\wedge \nu,dV_{\Sigma}\rangle dV_{\Sigma\cap \delta}}{det^{1/2}(\Pi^*_{\Sigma}\Pi_{\Sigma})} = \sum_a\epsilon_a\int_{(\Sigma\cap \delta)_a }tr\,Q,$$
where $\epsilon_a = \pm 1$ and $\Sigma\cap \delta = \cup_a (\Sigma\cap \delta)_a$ is the decomposition into disjoint components.

\subsection{Manifolds with boundary}
Because we have an analytic proof of the signature theorem which does not require the linear algebra of the index theorem, it is natural to consider the signature index form for manifolds with boundary with {\em local} boundary conditions. The well known theorems on topological obstructions to local boundary conditions for the signature operator \cite{AB} do not directly apply to our calculations, but one cannot hope to avoid the appearance of the $\eta$ invariant. In this subsection we see how this spectral term emerges from local boundary conditions. 

Let $M^{4k}$ be a compact manifold with boundary $Y$. Assume that near $Y$, $M$ is isometric to a product $[0,l]\times Y$. Then we may identify the integral of the index form over $\delta$ with the integral over $M$ of 
$$\tilde s_t:= -\sum_{i} h_{i}^{2k}(x)\wedge  h_{i}^{2k}(x)
- \sum_{\lambda}e^{-t\lambda}\lambda^{-1}\ast d\ast \psi_{\lambda}^{2k+1}(x)\wedge \ast d\ast\psi_{\lambda}^{2k+1}(x)$$
$$
- \sum_{\lambda}e^{-t\lambda}\lambda^{-1} d\psi_{\lambda}^{2k-1}(x)\wedge d\psi_{\lambda}^{2k-1}(x).$$
This eigen decomposition is not well defined until we fix boundary conditions. For simplicity we impose Dirichlet boundary conditions.  Then 
$$-\int_M\tilde s_t = \int_M \sum_{i} h_{i}^{2k}(x)\wedge  h_{i}^{2k}(x)
+ \sum_{\lambda}e^{-t\lambda}\lambda^{-1}\int_Y \ast \psi_{\lambda}^{2k+1}(x)\wedge d\ast\psi_{\lambda}^{2k+1}(x)$$
$$
+ \sum_{\lambda}e^{-t\lambda}\lambda^{-1} \int_Y\psi_{\lambda}^{2k-1}(x)\wedge d\psi_{\lambda}^{2k-1}(x).$$ The Dirichlet boundary conditions imply $\psi_{\lambda}^{2k-1}$ pulled back to $Y$ vanishes. Hence we are left with 
$$-\int_M\tilde s_t = \int_M \sum_{i} h_{i}^{2k}(x)\wedge  h_{i}^{2k}(x)
+ \sum_{\lambda}e^{-t\lambda}\lambda^{-1}\int_Y \ast \psi_{\lambda}^{2k+1}(x)\wedge d\ast\psi_{\lambda}^{2k+1}(x).$$
From the usual heat equation asymptotics we have 
\begin{equation}\label{bndry}\int_ML_k(TM) = \int_M \sum_{i} h_{i}^{2k}(x)\wedge  h_{i}^{2k}(x)
+ lim_{t\rightarrow 0}\sum_{\lambda}e^{-t\lambda}\lambda^{-1}\int_Y \ast \psi_{\lambda}^{2k+1}(x)\wedge d\ast\psi_{\lambda}^{2k+1}(x),
\end{equation}
where $L_k$ denotes the Hirzebruch $L$ polynomial. 
Thus we see that even with local boundary conditions, we obtain a spectral term similar to the eta invariant, although it is global on M rather than global on $Y$. Observe that because $\ast$ does not preserve Dirichlet boundary conditions (in fact it exchanges them with Neumann), we cannot diagonalize $\ast$ in the space of harmonic forms satisfying Dirichlet boundary conditions. Hence, the harmonic term in (\ref{bndry}) is not the signature and is not metric independent. On the other hand, if $H^{2k}(Y) = H^{2k-1}(Y) = 0,$ then it is easy to prove that 
$$\lim_{l\rightarrow\infty}\int_M \sum_{i} h_{i}^{2k}(x)\wedge  h_{i}^{2k}(x) = sig(M).$$
Hence, in the limit as the length $l$ of the collar tends to $\infty$, the  spectral term will converge to the eta invariant, by the APS signature theorem \cite{APS}. (The convergence of the spectral term can be proved without invoking the APS theorem, but the computations are effectively the same as in the proof of that theorem.)
\section{Spinors}
The preceding section suggests that the fixed point theorem for spinors should also extend to conformal correspondences. We now briefly examine such an extension. Let $M^{2m}$ be a compact, oriented Riemannian spin manifold. Let $S$ denote the bundle of spinors on $M$. Let $V\subset M\times M$ be an $n-$ dimensional smooth submanifold. Let 
$(x,y)\in V$ satisfy $d\pi_1 T_{(x,y)}V = T_xM$ and $d\pi_2 T_{(x,y)}V = T_yM$.  Then the projections determine a map
$$A(x,y):=d\pi_2\circ d\pi_1^{-1}:T_xM\rightarrow T_yM.$$
Let $P_{SO}(M)$ denote the oriented frame bundle of $M$, and let $P_{Spin}(M)$ denote the principal spin bundle. When $A=\mu(x,y) k(x,y)$, with $k(x,y)$ an isometry and $\mu(x,y) >0$ a scalar, then $k$ also induces a map from $(P_{SO}M)_x\rightarrow (P_{SO}M)_y$.  We call $V$ a {\em spin conformal correspondence} if it is a conformal correspondence, and there exists a continuous choice of lift of the action of $k(x,y)$ from the frame bundles to a map $k_{spin}:(P_{Spin}M)_x\rightarrow (P_{Spin}M)_y$, $(x,y)\in V$. A map between principal spin bundles induces a corresponding map of spin bundles 
$k_S:S_x\rightarrow S_y$.  Taking into account the entire conformal map $A$, we define 
$A_S :S_x\rightarrow S_y$ to be (see \cite[p.133]{LM})
$$A_S = \mu^{(1-2m)/2}k_S.$$
Then letting $D$ denote the Dirac operator on spinors (see \cite[Theorem 5.24]{LM}),
$$A_S D = D A_S.$$
With these preliminaries, we can extend the Lefschetz fixed point theorem to spin conformal correspondences. Let $\{\psi_{\lambda}\}_{\lambda}$ be an orthonormal basis of eigenspinors with $D^2$ eigenvalue $\lambda$.  Let $\{h_i\}_i$ be an orthonormal basis of harmonic spinors. Let $\tau_S$ now denote Clifford multiplication by the volume form, acting on spinors.  The index form for the spinor Dirac operator on a conformal correspondence $V$ is  
$$\sigma_t = \sum_{\lambda} e^{-t\lambda}\langle A_S(x,y)\psi_{\lambda}(x),\tau \psi_{\lambda}(y)\rangle dvol_M(y).$$
We expand this as 
$$\sigma_t = \sum_i\langle A_S(x,y)h_i(x),h_i(y)\rangle dvol_M(y) +  \sum_{\lambda\not = 0} \mu^{-1}e^{-t\lambda}\langle A_S(x,y)D\psi_{\lambda}(x),\tau D\psi_{\lambda}(y)\rangle dvol_M(y).$$
At the origin of a geodesic normal coordinate system, we write the $t-$ dependent term as 
$$\sum_{\lambda\not = 0} \lambda^{-1}e^{-t\lambda}\frac{\partial}{\partial y^j}\langle A_S(x,y)D\psi_{\lambda}(x),\tau c(dy^j)\psi_{\lambda}(y)\rangle dvol_M(y)$$
$$- \sum_{\lambda\not = 0} \lambda^{-1}e^{-t\lambda}\langle \nabla _{\frac{\partial}{\partial y^j}}A_S(x,y)D\psi_{\mu}(x),\tau c(dy^j)\psi_{\lambda}(y)\rangle dvol_M(y) $$
$$ = \sum_{\lambda\not = 0} \lambda^{-1}e^{-t\lambda}\frac{\partial}{\partial y^j}\langle A_S(x,y)D\psi_{\lambda}(x),\tau c(dy^j)\psi_{\lambda}(y)\rangle dvol_M(y)$$
$$- \sum_{\lambda\not = 0} \lambda^{-1}e^{-t\lambda}\langle DA_S(x,y)D\psi_{\lambda}(x),\tau \psi_{\lambda}(y)\rangle dvol_M(y) $$
$$ = \sum_{\lambda\not = 0} \lambda^{-1}e^{-t\lambda}\frac{\partial}{\partial y^j}\langle A_S(x,y)D\psi_{\lambda}(x),\tau c(dy^j)\psi_{\lambda}(y)\rangle dvol_M(y)$$
$$- \sum_{\lambda\not = 0} e^{-t\lambda}\langle A_S(x,y)\psi_{\lambda}(x),\tau \psi_{\lambda}(y)\rangle dvol_M(y).$$
Thus we see that the $t-$dependent term can be written 
$$\sum_{\lambda\not = 0} \frac{e^{-t\lambda}}{2\lambda}d[\langle A_S(x,y)D\psi_{\lambda}(x),\tau c(dy^j)\psi_{\lambda}(y)\rangle i_{\frac{\partial}{\partial y^j}}dvol_M(y)].$$
Here $i_X$ denotes interior multiplication by $X$.  Thus we see that the $t-$dependent term is once again exact.  We now obtain the following proposition.
\begin{proposition}
Let $M^{2m}$ be a compact, oriented Riemannian spin manifold.  Let $V$ be a spin conformal correspondence intersecting $\delta$ transversely. Then 
$$\int_V\sum_i\langle A_S(x,y)h_i(x),h_i(y)\rangle dvol_M(y) = \sum_{(b,b)\in V\cap \delta}\nu_{spin}(T_{(b,b)}V),$$
where for $V$ locally the graph of a conformal map $f = \mu k$ near $(b,b)$,
$$\nu_{spin}(T_{(b,b)}V) = \pm i^{m}2^{-m}\frac{\mu^{\frac{1}{2}-m}|det(I-k)|}{|det(I-\mu k)|}\prod_{j=1}^m cosec(\theta_j/2).$$
\end{proposition}
\section{Lagrangian geometry}
Suppose now that that $M$ is symplectic with symplectic form $\omega$. Suppose that $J$ is a compatible almost complex structure so that $\omega(JX,JY) = \omega(X,Y)$, for any vectorfields $X,Y$ on $M$, and $\omega(\cdot,J\cdot)$ defines a metric. The triple $(M,\omega,J)$ is called an {\em almost Kahler} manifold. Let $V$ be a coistropic submanifold of $M$. Recall that this means that for every $p\in V$,  the annihilator $A(T_pV)$ of  $T_pV$, defined by 
$$A(T_pV) =  \{v\in T_pM: \omega(v,\cdot)_{|T_pV} = 0\}$$
 is a subspace of $T_pV$. 
Let $N$ denote the normal bundle of $V$. Let $e_1,\cdots, e_l$ be an orthonormal basis of $A(T_pV)$. Then $Je_1,\cdots, Je_l$ span $N_p$ since they are contained in $N_p$ and $dim A(T_pV) + dim T_pV = dim T_p(M\times M).$ Consequently, we can complete the $e_j$ to an oriented orthonormal basis of $T_pV$, adding 
vectorfields $e_{l+1},Je_{l+1},\cdots, e_{l+p},Je_{l+p}$. Let 
$w^1,\cdots, w^l,w^{l+1},Jw^{l+1}\cdots,w^{l+p},Jw^{l+p}$ be the dual coframe. Consider a form $f:=w^I\wedge Jw^K.$
Observe that $\ast f_{|V} = 0$, unless $I = \emptyset$ and $K = \{1,\cdots,l\}$. In the latter case, we observe that 
$$\ast (Jw^1\wedge \cdots \wedge Jw^l) = (-1)^{l(l+1)/2}w^1 \wedge\cdots\wedge w^l\wedge w^{l+1}\wedge Jw^{l+1}\wedge\cdots \wedge w^{l+p}\wedge Jw^{l+p}.$$
We rewrite this as 
$$\ast f = \frac{1}{p!}L^p(-1)^{l(l-1)/2}Jf,$$
where $L$ denotes exterior multiplication by $\omega$. 
In general, we find that restricted to $V$,
$$\ast = \frac{1}{p!}L^p(-1)^{l(l-1)/2}J.$$
This suggests considering a modification of the signature index density, replacing 
$$s_t  = \sum_q(-1)^{q(n-q)}\sum_{\lambda}\tau \phi_{\lambda}^q(x)\wedge \ast\phi_{\lambda}^q(y)$$
by 
$$Js_t  = \sum_q(-1)^{q(n-q)}\sum_{\lambda}J\tau \phi_{\lambda}^q(x)\wedge \ast\phi_{\lambda}^q(y).$$
Then the $t$ dependent summands in the index density arise from exact $\phi$ and coexact $\phi$. These contribute terms of the form (up to constant factors) 
$\sum_{\lambda}J\tau \lambda^{-1}e^{-t\lambda}db_{\lambda}(x)\wedge \ast db_{\lambda}(y)$
and 
$\sum_{\lambda}J \lambda^{-1}e^{-t\lambda}d\beta_{\lambda}\wedge d\beta_{\lambda}.$
 Up to additional constant factors, for $V_1$ coisotropic, these restrict to $V_1\times V_2\subset M\times M$ as 
$\sum_{\lambda}\lambda^{-1}e^{-t\lambda}L^p db_{\lambda}(x)\wedge \ast db_{\lambda}(y)$
and 
$\sum_{\lambda}\lambda^{-1}e^{-t\lambda}J d\beta_{\lambda}(x)\wedge d\beta_{\lambda}(y).$ These integrate to $0$ for $V_i$ closed.  The $t-$invariant summand of $Js_t$ is therefore
$\sum_p(-1)^{p(n-p)}\sum_{i}J\tau h_i^p(x)\wedge \ast h_i^p(y).$

Applying (\ref{append}) yet again gives 
\begin{proposition} Let $V_1$ and $V_2$ be 2 transversely intersecting $m+q-$dimensional compact submanifolds of a compact almost Kahler manifold $M^{2m}$, with $V_1$ coisotropic. Let $\{h_j^q\}_{j=1}^{b_q}$ be an orthonormal basis of harmonic $q-$ forms of $M$. Then
\begin{equation}
\sum_j\int_{V_1}Jh_j^q\int_{V_2} h_j^q = \int_{V_1\cap V_2}\frac{<dvol_{V_1},Jdvol_{V_2}>}{det^{1/2}(2\Pi^*_{V_1\times V_2,b}\Pi_{V_1\times V_2,b})}dV_{V_1\cap V_2}(b).
\end{equation}

\end{proposition}

\begin{remark}This proposition can also be deduced using the Gauss-Bonnet index form. To see this write 
$$
\frac{-(-1)^{(m-q)(m-q-1)/2}}{q!}\sum_j\int_{V_1}Jh_j^q\int_{V_2} h_j^q = (-1)^{m-q}\sum_j\int_{V_1}L^q h_j^{m-q}\int_{V_2} \ast h_j^{m-q}$$
$$= \int_{V_1\times V_2}\omega^q \wedge \sum_p(-1)^{p}e_t^p.$$
This can be computed using Theorem \ref{fatgb}. Using the fact that $V_1$ is coisotropic we recover the equality of the proposition.   
\end{remark}
\section{Special Lagrangians}
Let $M^{2m}$ be a compact Calabi - Yau manifold. 
Let  $V = V_1\times V_2$, with each $V_i$ a special Lagrangian manifold. Recall $M$ Calabi-Yau implies that there is a nonvanishing holomorphic $(m,0)$ form, $\Omega$ on $M$. $V_i$ special Lagrangian means that $V_i$ is Lagrangian and $\Omega|_{V_j} = e^{i\phi_j}dvol_{V_j}$, for some constant phase $e^{i\phi_j}$.  It is obvious that the $t-$independence of the Riemann-Roch index form must fail for $V$,  as the integral reduces to 
$$(-1)^m\int_{V_1}h^{0,m}\int_{V_2}\bar h^{0,m} + (-1)^m\sum_{\lambda}e^{-t\lambda}\lambda^{-1}\int_{V_1}\db b_{\lambda}^{0,m-1}\int_{V_2}\partial \bar b_{\lambda}^{0,m-1}.$$

In the $t\rightarrow \infty$ limit, this reduces to 
$$(-1)^me^{-i(\phi_1-\phi_2)}vol(V_1)vol(V_2)vol(M)^{-1}.$$

On the other hand, the $t\rightarrow 0$ limit (\ref{nRR}) yields the expression,

$$lim_{t\rightarrow 0}\int_{V_1\times V_2}\sum_q e_t^{0,q} = (-1)^m\sum_{b\in V_1\cap V_2}\frac{e^{-i(\phi_1-\phi_2)}}{det^{1/2}(2\Pi^*_{V_1\times V_2,b}\Pi_{V_1\times V_2,b})}.$$

It seems unlikely that these two limits should coincide. We record their relation:
\begin{proposition}Let $M$ be a compact Calabi Yau. Let $V_1$ and $V_2$ be two special Lagrangian submanifolds of $M$ intersecting transversely. Then 
$$e^{i(\phi_1-\phi_2)}vol(V_1)vol(V_2)vol(M)^{-1} - \sum_{b\in V_1\cap V_2}\frac{e^{-i(\phi_1-\phi_2)}}{det^{1/2}(2\Pi^*_{V_1\times V_2,b}\Pi_{V_1\times V_2,b})}$$
$$= \sum_{\lambda}\lambda^{-2}\int_{V_1}\db b_{\lambda}^{0,m-1}\int_{V_2}\partial \bar b_{\lambda}^{0,m-1}.$$
\end{proposition}
If we weaken Calabi-Yau to Kahler and special Lagrangian to Lagrangian, then this relation generalizes to
$$\sum_i\int_{V_1}h_i^{0,m}\int_{V_2}\bar h_i^{0,m} - \sum_{b\in V_1\cap V_2}\nu_{RR}(T_{(b,b)}(V_1\times V_2))$$
$$= \sum_{\lambda}\lambda^{-2}\int_{V_1}\db b_{\lambda}^{0,m-1}\int_{V_2}\partial \bar b_{\lambda}^{0,m-1}.$$  
On an abelian variety $M=C^m/L$ , the spectral term has the simple form 
$$\frac{1}{vol(M)}\sum_{k\in L^*\setminus 0}|k|^{-2}\int_{V_1}e^{ik\cdot x}dv_x\int_{V_2}e^{-ik\cdot y}dv_y,$$
and invites data mining. For example, performing Hecke type operations on the lattice leads to many relations. Even simpler is to consider a $2m$ parameter family of translations $V_{2s} = V_2 + s$. Then integrating the equality in $s$ gives 
$$e^{i(\phi_1-\phi_2)}vol(V_1)vol(V_2)vol(M)^{-1}\int_S ds $$
$$ - \int_S\sum_{b(s)\in V_1\cap V_{2s}}\frac{e^{-i(\phi_1-\phi_2)}}{det^{1/2}(2\Pi^*_{V_1\times V_{2s},b}\Pi_{V_1\times V_{2s},b})}ds$$
$$= \sum_{k\in L^*\setminus 0}\frac{\int_Se^{-ik\cdot s}ds}{vol(M)|k|^2}\int_{V_1}e^{ik\cdot x}dv_x\int_{V_2}e^{-ik\cdot y}dv_y.$$
Choosing $S$ to be $C^m/L$, annihilates the spectral term, leaving 
\begin{equation}\label{average}
e^{i(\phi_1-\phi_2)}vol(V_1)vol(V_2) = \int_M\sum_{b(s)\in V_1\cap V_{2s}}\frac{e^{-i(\phi_1-\phi_2)}}{det^{1/2}(2\Pi^*_{V_1\times V_{2s},b}\Pi_{V_1\times V_{2s},b})}ds.\end{equation}

Similar computations can be performed on other homogeneous spaces and for other index forms. For example, if $M = G/K$ is a compact hermitian symmetric space and if $V_1$ and $V_2$ are Lagrangian submanifolds, we can use the group $G$ to translate $V_2$. Then the integral of the shifted spectral term over $G$ becomes 
$$\sum_{\lambda}\lambda^{-2}\int_Gdg\int_{V_1}\db b_{\lambda}^{0,m-1}\int_{gV_2}\partial \bar b_{\lambda}^{0,m-1} = \sum_{\lambda}\lambda^{-2}\int_{V_1}\db b_{\lambda}^{0,m-1}\int_{V_2}\partial \bar (\int_Gdg g^*b_{\lambda}^{0,m-1}).$$  
The form $\int_Gdg g^*b_{\lambda}^{0,m-1}$ is invariant and therefore vanishes as it must be harmonic but $\lambda\not = 0.$
So we obtain 
\begin{theorem}Let $M$ be a compact Hermitian symmetric space. Let $V_1$ and $V_2$ be two Lagrangian submanifolds of $M$ intersecting transversely. Then 
$$\sum_i\int_{V_1}h_i^{0,m}\int_{V_2}\bar h_i^{0,m} = Vol(G)^{-1}\int_Gdg\sum_{b\in V_1\cap gV_2}\nu_{RR}(T_{(b,b)}(V_1\times gV_2)).$$  
\end{theorem}
The same argument gives 

\begin{theorem}Let $M= G/K$ be a compact symmetric space. Let $V_1$ and $V_2$ be two middle dimensional submanifolds of $M$ intersecting transversely. Let $\{h_i\}_i$ be an orthonormal basis for the harmonic forms in middle dimension. Then 
$$\sum_i\int_{V_1}h_i \int_{V_2} h_i  = Vol(G)^{-1}\int_Gdg\sum_{b\in V_1\cap gV_2}\nu_{sig}(T_{(b,b)}(V_1\times gV_2)).$$  

\end{theorem}

Thus the average of the local signature numbers gives the intersection number of $V_1$ and the poincare dual of $V_2$. 
\section{Appendix: heat equation asymptotics}

Let $E$ be a Dirac bundle over a smooth compact Riemannian manifold $M^n$. (See \cite[p.114]{LM}). 
Let $D$ be the associated generalized Dirac operator. We review here the construction of an approximation to $e^{-tD^2}$. 
The Cauchy integral formula
$$e^{-tD^2} = \frac{-1}{2\pi i}\int_{\gamma}e^{-\lambda}(tD^2-\lambda)^{-1}d\lambda$$
reduces the construction to approximating $(tD^2-\lambda)^{-1}$. 
The standard method of approximation (see \cite{Gil} or \cite{BGV}), which we will follow here, is to construct an approximation in coordinate neighborhoods. These local approximations are then patched together using partitions of unity and auxillary cutoff functions. We will suppress this latter patching step in our discussion.  

 Fix $y\in M$ and geodesic coordinates centered at $y$. Define
$$P_{\lambda,N}f(x) = \int e^{2\pi i(x-y)\cdot u}\sum_{j=0}^N(4\pi^2t|u|^2-\lambda)^{-j-1}a_j(x,y)f(y)dydu,$$
with $a_0=Id$ in our choice of local frames. 
The remaining $a_j$ are chosen inductively with  
$$(4\pi^2t|u|^2-\lambda)^{-j}a_j(x,y) = -(tD_x^2-4\pi itu^k\nabla_k)(4\pi^2t|u|^2-\lambda)^{-j}a_{j-1}(x,y),$$
for $1\leq j\leq N$. This gives the recipe 
$$(4\pi^2t|u|^2-\lambda)^{-j-1}a_j(x,y) = (-t)^j(4\pi^2t|u|^2-\lambda)^{-1}[(D_x^2-4\pi iu^k\nabla_k)(4\pi^2t|u|^2-\lambda)^{-1}]^ja_{0}(x,y).$$
With this choice, 
$$(tD^2-\lambda)P_{\lambda,N}f(x) = \int (tD_x^2-\lambda)e^{2\pi i(x-y)\cdot u}\sum_{j=0}^N(4\pi^2t|u|^2-\lambda)^{-j-1}a_j(x,y)f(y)dydu$$
$$= \int e^{2\pi i(x-y)\cdot u}(tD_x^2-4\pi itu^k\nabla_k + 4\pi^2|u|^2-\lambda)\sum_{j=0}^N(4\pi^2t|u|^2-\lambda)^{-j-1}a_j(x,y)f(y)dydu$$

$$= f(x) + \int e^{2\pi i(x-y)\cdot u}[(tD_x^2-4\pi itu^k\nabla_k)(4\pi^2t|u|^2-\lambda)^{-N-1}a_N(x,y)]f(y)dydu,$$

Inserting this back into our expression for $e^{-tD^2}$ gives, for a suitable curve $\gamma$ in $\IC$ surrounding the real axis, 
the approximate heat kernel 
$$p^N_t(x,y) = \frac{-1}{2\pi i}\int_{\gamma}e^{-\lambda}\int e^{2\pi i(x-y)\cdot u}\sum_{j=0}^N(-t)^j(4\pi^2t|u|^2-\lambda)^{-1}[(D_x^2-4\pi iu^k\nabla_k)(4\pi^2t|u|^2-\lambda)^{-1}]^ja_{0}(x,y)dud\lambda $$
$$ = \frac{-t^{-n/2}}{2\pi i}\int_{\gamma}e^{-\lambda}\int e^{2\pi i(x-y)\cdot u/\sqrt{t}}\sum_{j=0}^N(-t)^j(4\pi^2|u|^2-\lambda)^{-1}[(D_x^2-4\pi it^{-1/2}u^k\nabla_k)(4\pi^2|u|^2-\lambda)^{-1}]^ja_{0}(x,y)dud\lambda.$$
The error term $p_t-p_t^N$ has trace class norm which is decreasing faster than $O(t^{N/4})$, (not sharp) for $N$ large and $t\rightarrow 0$. 
Expand
$$[(D_x^2-4\pi it^{-1/2}u^k\nabla_k)(4\pi^2|u|^2-\lambda)^{-1}]^ja_0 = \sum_{l,J,p}(4\pi^2|u|^2-\lambda)^{-l}u^Jt^{-p/2}a_{j,l,J,p}(x,y).$$
In this expansion, we note that $p\leq j$.
Inserting this into our expression for $p^N_t(x,y)$, changing the order of integration, and performing the contour integral gives 

$$p^N_t(x,y) = \int e^{-4\pi^2|u|^2}e^{2\pi i(x-y)\cdot u/\sqrt{t}}\sum_{j=0}^N\sum_{l,J,p}\frac{(-1)^ju^Jt^{j-p/2-n/2}}{l!}a_{j,l,J,p}(x,y)du$$
$$ = \sum_{j=0}^N\sum_{l,J,p}[(\frac{\partial}{2\pi i\partial x})^J  e^{-|x-y|^2/4t}](4\pi)^{-n/2}\frac{(-1)^jt^{j+|J|/2-p/2-n/2}}{l!}a_{j,l,J,p}(x,y).$$
Let $\tau_E$ denote an involution of $E$, preserving the fibers. We wish to study integrals of $\tau_E p_t$ over submanifolds of $M\times M$, but in general $\tau_E p_t$ is a section of $\pi_1^*E\otimes\pi_2^*E^*.$ Hence we must first fix a section 
$q\in \Gamma(Hom (\pi_1^*E\otimes\pi_2^*E^*,\bigwedge^{\cdot}T^*(M\times M)))$ and consider instead integrals of 
$q(\tau_E p_t^N).$ Let $W$ be a smooth submanifold of $M\times M$ intersecting $\delta$ transversely, in a submanifold $S$. 
Then Gaussian decay gives 
$$lim_{t\rightarrow 0}\int_W q(\tau_E p_t) = lim_{t\rightarrow 0}\int_W q(\tau_E p_t^N) = lim_{t\rightarrow 0}\int_{N_{\epsilon}(S)} q(\tau_E p_t^N),$$
for any tubular neighborhood $N_{\epsilon}(S)$ of $S$ in $W$. In a neighborhood $V$ of $(b,b)\in S$, pick a local orthonormal frame $\{e_j\}_j$ for the normal bundle $N_S$ of $S$ in $W$. Given a coordinate map $X$ for $V$, we obtain coordinates for the tubular neighborhood of $V$ via the map  
$$(v,s)\rightarrow exp_{X(v)}(s^je_j).$$ 
The Gaussian $e^{-|x-y|^2/4t}$ can be written $e^{-d(y,x)^2/4t}$, since the exponential map preserves radial distances. Observe that 
$$d(\pi_1(exp_{X(v)}(s^je_j)), \pi_2(exp_{X(v)}(s^je_j)))^2 = 2|\Pi_W(s^je_j)|^2 + O(|s|^4),$$
where we recall that $\Pi_W$ denotes the projection of the normal bundle of $S$ in $W$  onto the normal bundle of $\delta$. 
This follows from the fact that the result is true in Euclidean space and the exponential map distorts distances by at most $O(|s|^2)$. 
Let $\phi(b,s)$ denote the $O(|s|^4)$ correction term. More generally, we can write 
$(x-y) = 2^{1/2}\Pi_W(s) + \phi_2(b,s),$ with $\phi_2(b,s)\in O(|s|^2),$ and 
$dV_W = dV_S\wedge (1+\phi_3(b,s))ds,$ with $\phi_3(b,s)\in O(|s|),$
 Write $(\frac{\partial}{2\pi i\partial x})^J  e^{-|x-y|^2/4t} = w_J(x-y,t)e^{-|x-y|^2/4t},$ for some polynomial $w_J = \sum_{2b-|A| = |J|}w_{J,A,b}(x-y)^At^{-b}$. Then we wish to compute 

$$lim_{t\rightarrow 0}\int_{V\times B_{\epsilon}(0)}\sum_{j=0}^N\sum_{l,J,p}w_J(2^{1/2}\Pi_W(s)+\phi_2(b,s),t)  e^{-(2|\Pi_W s|^2+\phi(b,s))/4t}$$
$$\times \frac{(-1)^jt^{j+|J|/2-p/2-n/2}}{(4\pi)^{n/2}l!}\langle q(\tau_E a_{j,l,J,p}(b,s)),dV_W\rangle dV_S(1+\phi_3(s,b))ds $$
$$ = lim_{t\rightarrow 0}\int_{V\times B_{t^{-1/2}\epsilon}(0)}\sum_{j=0}^N\sum_{l,J,p}w_J(2^{1/2}\Pi_W(s)+t^{-1/2}\phi_2(b,t^{1/2}s),1)  e^{-(2|\Pi_W s|^2+\phi(b,t^{1/2}s)/4t)}$$
$$\times \frac{(-1)^jt^{j-p/2-n/2+c/2}}{(4\pi)^{n/2}l!}\langle q(\tau_E  a_{j,l,J,p}(b,t^{1/2}s)),dV_W\rangle dV_S(1+\phi_3(b,t^{1/2}s))ds,$$
where $c$ denotes the codimension of $S$ in $W$. 

In this paper, we have only treated the case where $c=n$. In this case, recalling that $p\leq j$, the preceding vanishes unless $j=p=|J|=0=l$ and then reduces to 
$$\int_{V\times \IR^n} e^{-2|\Pi_W s|^2}\frac{1}{(4\pi)^{n/2}}\langle q(\tau_E),dV_W\rangle dV_Sds = \int_{V}\frac{\langle q(\tau_E),dV_W\rangle }{det^{1/2}(2\Pi_W^*\Pi_W(b))}dV_S.$$
We record this as a proposition.
\begin{proposition}\label{append}Let $E$ be a Dirac bundle over a smooth compact Riemannian manifold $M^n$. Let $W$ be a smooth $n+p$ dimensional submanifold of $M\times M$ which intersects $\delta$ transversely. Let $q\in \Gamma(Hom(\pi_1^*E\otimes\pi_2^*E^*,\bigwedge^{\cdot}T^*(M\times M)))$. Let $\tau_E$ be an automorphism of $E$ which is an involution. Then
\begin{equation}lim_{t\rightarrow 0}\int_W q(\tau_E p_t) = \int_{W\cap\delta}\frac{\langle q(\tau_E),dV_W\rangle }{det^{1/2}(2\Pi_W^*\Pi_W(b))}dV_{W\cap\delta}.
\end{equation}
\end{proposition}

{\bf Example 1}: As an example, suppose that $E$ is $\bigwedge ^{\cdot}T^*M$, and identify $E^*$ with differential forms via the Hodge star operator. Let $q=q_1$ be the identity map with respect to these identifications. Then if $E=\bigwedge ^{\cdot}T^*M$ and if $\{\omega_j\}_j$ is an oriented orthonormal basis of $T_b^*M$, 
$q(Id)(b,b) = \sum_J\pi_1^*\omega^{J}\wedge\pi_2^*\ast\omega^J$.  Suppose that $\tau_E = (-1)^{deg }$. 
Then $$q(\tau_E) = \sum_J(-1)^{|J|}\pi_1^*\omega^{J}\wedge\pi_2^*\ast\omega^J 
= (\pi_1^*\omega^1-\pi_2^*\omega^1)\wedge\cdots\wedge (\pi_1^*\omega^n-\pi_2^*\omega^n).$$
Let $\nu$ denote the volume form of the fiber of normal bundle to $\delta$, determined by the metric. 
Then 
$$q(\tau_E) = 2^{n/2}\nu.$$
Suppose, as in the classical Lefschetz theorem that $W$ is the graph of $f$, where $f$ has only isolated nondegenerate fixed points. Then if $(b,b)\in W\cap\delta,$
$T_{(b,b)}W$ can be identified with the graph of $df_b$. Let $\{v_j\}_j$ be an oriented orthonormal eigenbasis for $df_b^*df_b$, with $df_b^*df_bv_j = \lambda_jv_j$. Then 
$\{(1+\lambda_j)^{-1/2}(df_bv_j,v_j)\}_j$ is an orthonormal basis of $T_{(b,b)}W$. Hence
$\langle q(\tau_E),dV_W\rangle $ can be computed by evaluating $2^{n/2}\nu$ on this basis, yielding 
$$\langle q(\tau_E),dV_W\rangle  = det(I-df_b)\prod_j(1+\lambda_j)^{-1/2}.$$
On the other hand, we compute 
$$\langle 2\Pi_W^*\Pi_W(1+\lambda_j)^{-1/2}(df_bv_j,v_j),(1+\lambda_k)^{-1/2}(df_bv_k,v_k)\rangle = 
(1+\lambda_j)^{-1/2}(1+\lambda_k)^{-1/2}\langle v_j-df_bv_j,v_k-df_bv_k\rangle $$
$$
 = (1+\lambda_j)^{-1/2}(1+\lambda_k)^{-1/2}\langle (I-df_b^*)(I-df_b)v_j,v_k\rangle.$$
So $det^{1/2}(2\Pi_W^*\Pi_W) = |det(I-df_b)|\prod_j(1+\lambda_j)^{-1/2}. $ Combining these gives the expected
\begin{equation}\label{example}
\frac{\langle q(\tau_E),dV_W\rangle }{det^{1/2}(2\Pi_W^*\Pi_W(b))} = sign\,det(I-df_b),
\end{equation}
and we see that $\frac{\langle q(\tau_E),dV_W\rangle}{det^{1/2}(2\Pi_W^*\Pi_W(b))}$ encodes the usual Lefschetz data. 

{\bf Example 2}: Let $E$ and $\tau_E$ be as in the preceding example, but choose $q_2\in \Gamma(Hom (\pi_1^*E\otimes\pi_2^*E^*,\bigwedge^{\cdot}T^*(M\times M)))$ to be 
$z\wedge q_{1}$, for some $p-$form $z$. Choose $W$ to be an $n+p$ dimensional submanifold of $M\times M$ intersecting $\delta$ transversely in a submanifold $S$. On $S$, the metric gives us an orthogonal factorization $dV_{W} = \nu_S\wedge dV_S$. Then 
$$\frac{\langle q_2(\tau_E),dV_W\rangle}{det^{1/2}(2\Pi_W^*\Pi_W(b))} = \frac{\langle z\wedge 2^{n/2}\nu,\nu_S\wedge dV_S\rangle}{det^{1/2}(2\Pi_W^*\Pi_W(b))} = \frac{\langle z,dV_S\rangle\langle2^{n/2}\nu,\nu_S\rangle}{det^{1/2}(2\Pi_W^*\Pi_W(b))}.$$
We have seen in the previous example that 
$$\frac{(2^{n/2}\nu,\nu_S)}{det^{1/2}(2\Pi_W^*\Pi_W(b))} = \pm 1.$$
Hence 
$$\frac{\langle q_2(\tau_E),dV_W\rangle}{det^{1/2}(2\Pi_W^*\Pi_W)}dV_S = \pm z_{|S},$$
where the sign is constant on connected components of $S$.  

{\bf Example 3}:  
In treating the fixed point theorem for the Dolbeault complex on $(0,q)$ forms, $q(\tau_E) = 2^{n/2}\nu$ is replaced by 
$q(\tau_E) = 2^{m/2}\nu_2,$ where $\nu_2$ is defined in \ref{nu2} and $m$ denotes the complex dimension of $M$. Let's consider the more general choice of $q = z\wedge 2^{m/2}\nu_2$, $z$ a closed $(p,p)$ form, arising when we consider currents of greater dimension than $M$. Then in the notation of section \ref{excesssub}, with $S = W\cap \delta$, 
$$\frac{\langle q(\tau_E),dV_W\rangle }{det^{1/2}(2\Pi_W^*\Pi_W(b))}dV_{S} 
= \frac{\langle z\wedge 2^{m/2}\nu_2,\nu_S^W\wedge dV_S \rangle }{det^{1/2}(2\Pi_W^*\Pi_W(b))}dV_{S}.$$
Suppose $W$ is given locally as the projection onto the two $M$ factors of the graph of a function 
$$F:M\times B\rightarrow M,$$
where $B\subset \pi_1(S)$ denotes a small ball. 
Thus locally $W$ has the form $\{(z,F(z,t)):z\in U\subset M, t\in B\}. $
By making a change of variables in $U\times B$ if necessary, we can assume that at a given fixed point $(s,s)\in S$, $F(s,0) = s$, $(dF,I)T_0B = T_{s,s}S,$ and $(dF_{(s,0)}u,u)\perp (dF_{(s,0)}v,v)$, for $u\in T_sU$ and $v\in T_0B$. Let $f(x) = F(x,0)$. Choose again a basis $\{v_j\}_j$ for $T_{s}M$ which is eigen for $df_{s}^*df_{s}$: $df_{s}^*df_{s}v_j = \lambda_jv_j$. 
Then the tangent space to $W$ at $(s,s)$ has orthonormal basis $\{(1+\lambda_j)^{-1/2}(dfv_j,v_j)\}_j$ union an orthonormal basis $\{u_a\}_a$ for $T_{(s,s)}S$. With these notations, 
$det^{1/2}(2\Pi_W^*\Pi_W)$ is now the same as $det^{1/2}(2\Pi_{\Gamma_f}^*\Pi_{\Gamma_f})$, 
where $ \Gamma_f $ denotes the graph of $f$. In particular, 
$$det^{1/2}(2\Pi_W^*\Pi_W) = |det(I-df)|\prod_j(1+\lambda_j)^{-1/2}.$$
To proceed further, we need to decompose according to Hodge type
$$\nu_S^W = \nu_S^W(0,n)\wedge \nu_S^W(n,0),$$
and 
$$2^{m/2}\nu_2 = 2^{m/2}\nu_2(0,n)\wedge \nu_2(n,0).$$
Then we have  

$$\frac{\langle z\wedge 2^{m/2}\nu_2,\nu_S^W\wedge dV_S \rangle }{det^{1/2}(2\Pi_W^*\Pi_W)} = 
\frac{\langle 2^{m/2}\nu_2(0,n),\nu_S^W(0,n)\rangle \langle z\wedge \nu_2(n,0),\nu_S^W(n,0)\wedge dV_S \rangle }{|det(I-df)|\prod_j(1+\lambda_j)^{-1/2}}.$$
Computing as in Example 1, we find 

$$\frac{\langle 2^{m/2}\nu_2(0,n),\nu_S^W(0,n)\rangle }{|det(I-df)|\prod_j(1+\lambda_j)^{-1/2}} = \frac{\prod_j(1+\lambda_j)^{1/2}}{det_{\IC}(I-df)}.$$
This reduces us to computing 

$$ \sum_a\langle z_{p,a,0,p-a}\wedge \nu_2(n,0),\nu_S^W(n,0)\wedge dV_S \rangle $$
$$ = 
(-i)^p\prod_{j=1}^n(1+\lambda_j)^{-1/2} \sum_a z_{p,a,0,p-a}\wedge \nu_2(n,0)((dfv_1,v_1),\cdots,(dfv_n,v_n), u_1,\bar u_1,\cdots,u_p,\bar u_p), $$
where we have ordered our eigenbases so that the $v_j, j\leq n$ and $u_a$ , $a\leq p$ are type $(1,0)$. 
The orthogonality conditions on $(dF,I)$ imply that $d\pi_1(T_{s,s}S)$ must be contained in the $-1$ eigenspace of $df_s$. Hence we may assume that we have chosen our basis so that $u_j = \frac{1}{\sqrt{2}}(v_j,v_j)$. Hence the preceding reduces to 

$$ \sum_a\langle z_{p,a,0,p-a}\wedge \nu_2(n,0),\nu_S^W(n,0)\wedge dV_S \rangle $$
$$ = i^p\prod_{j=p+1}^n(1+\lambda_j)^{-1/2} \sum_a z_{p,a,0,p-a}(u_1,\bar u_1,\cdots, u_p,\bar u_p ) $$
$$ = 2^{p/2}\prod_{j=p+1}^n(1+\lambda_j)^{-1/2} \langle\sum_a z_{p,a,0,p-a},dV_S\rangle.$$

Combining these expressions and recalling the notation $z_0 = \sum_a z_{p,a,0,p-a}$ gives 

\begin{equation}\label{hololast}\frac{\langle z\wedge 2^{m/2}\nu_2,\nu_S^W\wedge dV_S \rangle }{det^{1/2}(2\Pi_W^*\Pi_W(b))}dV_{S} = 
\frac{2^{p}\langle z_{0},dV_S\rangle dV_S}{det_{\IC}(I-df)} = 
2^p\frac{\langle 2^{m/2}\nu_2,\nu_S^W\rangle}{det^{1/2}(2\Pi_W^*\Pi_W)} (z_{0})_{|S}.
\end{equation}
%
 \newcommand{\etalchar}[1]{$^{#1}$}
\def\polhk#1{\setbox0=\hbox{#1}{\ooalign{\hidewidth
  \lower1.5ex\hbox{`}\hidewidth\crcr\unhbox0}}}
\providecommand{\bysame}{\leavevmode\hbox to3em{\hrulefill}\thinspace}

%

\end{document}